\documentclass[11pt]{article}

\usepackage{amsmath}
\usepackage{amssymb}
\usepackage{amsthm}

\def\implies{\Longrightarrow}

\def\g{\mathfrak{g}}
\def\H{\mathcal{H}}
\def\V{\mathcal{V}}
\def\S{\mathcal{S}}
\def\W{\mathcal{W}}
\def\X{\mathfrak{X}}
\def\R{\mathbb{R}}
\def\Q{{\sf Q}}

\def\vectorfields#1{\X(#1)}

\def\fpd#1#2{{\displaystyle\frac{\partial #1}{\partial #2}}}

\def\lie#1{{\mathcal{L}}_{#1}}
\def\vf#1{\frac{\partial}{\partial #1}}
\def\vlift#1{#1^{\scriptscriptstyle{\mathrm{V}}}}

\def\hlift#1{#1^{\scriptscriptstyle{\mathrm{H}}}}
\def\conn#1#2#3{\setbox1=\hbox{$\scriptstyle{#2}{#3}$}%
\setbox2=\hbox to\wd1{$\hfil\scriptstyle{#1}\hfil$}
\Gamma^{\!\box2}_{\!\box1}}

\def\ad{\mathop{\mathrm{ad}}\nolimits}
\def\id{\mathop{\mathrm{id}}}
\def\im{\mathop{\mathrm{im}}}
\def\ker{\mathop{\mathrm{ker}}}
\def\rank{\mathop{\mathrm{rank}}}

\def\onehalf{{\textstyle\frac12}}

\newtheorem{thm}{Theorem}

\newtheorem{prop}{Proposition}
\newtheorem{cor}{Corollary}

\begin{document}

\title{Involutive distributions and dynamical systems of second-order type}

\author{T.\ Mestdag and M.\ Crampin\\
Department of Mathematics, Ghent University\\
Krijgslaan 281, S9, B--9000 Gent, Belgium\\ tom.mestdag@ugent.be, crampin@btinternet.com}

\date{}

\maketitle

{\small {\bf Abstract.} We investigate the existence of coordinate
transformations which bring a given vector field on a manifold
equipped with an involutive distribution into the form of a
second-order differential equation field with parameters.  We define
associated connections and we give a
coordinate-independent criterion for determining whether the vector
field is of quadratic type.  Further, we investigate the underlying
global bundle structure of the manifold under consideration, induced
by the vector field and the involutive distribution.  \\[2mm]
}

\section{Introduction}

In its most general sense, a (smooth) dynamical system, from the
geometrical point of view, is simply a vector field on some manifold.
But many dynamical systems of interest in physics and engineering
applications are more specialized than that:\ they are of second-order
type.  By a dynamical system, or vector field, of second-order type,
or a second-order differential equation field, we mean a vector field
$\Gamma$ on the tangent bundle $\tau:T(Q)\to Q$ of some configuration
manifold $Q$ with the property that $\tau_{*y}\Gamma=y$ for all $y\in
T(Q)$, so that in terms of coordinates $(x^i,y^i)$ where the $y^i$ are
the canonical fibre coordinates corresponding to coordinates $x^i$ on
$Q$,
\[
\Gamma=y^i\vf{x^i}+\Gamma^i(x,y)\vf{y^i}.
\]
It is of interest therefore, to find criteria for determining whether
a given dynamical system, which may be represented in some arbitrary
coordinates, is actually of second-order type, in that coordinates may
be found with respect to which it takes the form above.  This is a
problem which has both a local aspect, just described, and a global
one, which includes the question of whether the manifold on which the
dynamical system resides is in fact a tangent bundle.

In the recent article \cite{RR}, Ricardo and Respondek deal with a
version of the problem in the context of control theory.  Assume $M$
to be an even dimensional manifold, not known to be the tangent
manifold of some configuration manifold.  Under which conditions does
a coordinate change exist such that a given nonlinear control system
\[
\dot z^\alpha = F^\alpha(z) + u_r G^\alpha_r(z), \qquad
\alpha = 1,\ldots, 2n,
\]
is transformed into the form of a so-called `mechanical control
system', meaning a dynamical system of the form
\[
{\ddot q}^a = \conn abc(q) {\dot q}^b {\dot q}^c
+ P^a_b(q){\dot q}^b +Q^a(q) + u_r g^a_r(z), \qquad a=1, \ldots, n?
\]  The solution of the problem in \cite{RR} is cast in
terms of a certain vector space $\V$, which is a subspace of the
infinite dimensional vector space of vector fields on $M$, with
dimension exactly the half of the dimension of the manifold, and
which, among other properties, contains the control forces $G_r$ and
satisfies $[\V,\V]=0$.

In this paper, we will address a more general problem.  First of all,
we will not assume that the dimension of the manifold is even.  This
is mainly motivated by the observation that even for a system of
time-dependent second-order differential equations the manifold on
which the dynamics is described is odd-dimensional:\ it is the first
jet manifold of a spacetime manifold, or event space, fibred
over the real numbers (see e.g.\ \cite{CMS,MaPa}).
Further, we will allow some of the new coordinates to play simply the
role of parameters.  That is to say, we will not require the number
$n$ of second-order equations to be exactly the half of the dimension
of the manifold $M$.

Next to extending the results of \cite{RR} to a broader class of
manifolds, we will also make some conceptual modifications.  In a
nutshell, the results of \cite{RR} claim that if $\V$ (which is a
vector space constructed from the given control forces $G_r$ and from
$F$) satisfies certain conditions the vector field $F=F^\alpha
\partial/\partial z^\alpha$ (the so-called drift vector field)
transforms into an appropriate coordinate form, and, as a side-effect,
so does also the controlled dynamical field $F+u_rG_r$.  We will take
the space $\V$ to be the primary given object of our study, and ignore
that it was constructed from some given control forces.  Consequently,
we shift the attention to specific coordinate expressions for $F$, and
leave the control system given by the vector field $F+u_rG_r$ out of
the picture all together.  A second deviation is that for us $\V$ will
not be a vector space of vector fields, but rather the distribution it
generates.

In Section~2 we investigate under what conditions a given vector field
$F$ can be transformed into the coordinate expression of a
second-order differential equation field with possible parameters, in
the presence of an arbitrary involutive distribution $\V$ (of
arbitrary dimension).  Our framework has the advantage that it leaves
open the possibility that the transformed dynamics become either
autonomous or time-dependent.  We show in Section~3 how one can
associate various
connections to $F$, and we argue that these connections provide a
coordinate-independent method to express that the dynamics of $F$ is
of quadratic type (or of mechanical type, in the sense as above) in a
yet unknown set of coordinates.

Working with a distribution $\V$ has the further advantage that it
brings an associated almost tangent structure (and almost jet
structure) to the foreground.  These geometric structures find their
equivalence in standard tangent bundle and jet bundle geometry, but
they went unnoticed in \cite{RR}.  Based on results in
\cite{CT,DeF,MP2}
we further address in Section~4 the global issues that arise in this
context, such as e.g.\ the affine fibre bundle structure of $M$ and
the relation of $F$ to second-order differential equation fields on a
certain tangent or jet manifold.

In the last section we illustrate the theory in the context of a
Lagrangian system with an Abelian symmetry group, where a second-order
differential equation field with multiple parameters naturally shows
up.

\section{Local coordinate transformations}

Let $M$ be a manifold of dimension $m$ and $\V$ an involutive
distribution on $M$ of dimension $n$ such that $2n\leq m$.

We must be a little careful here about the meaning of the term
distribution and related terms.  A distribution on $M$ is of course a
choice of subspace of $T_z(M)$ at each $z\in M$, of constant
dimension, depending smoothly on $z$ in the sense that it admits local
smooth bases.  There is a related, but distinct, concept, which we may
call a vector field system.  A vector field system $\S$ on $M$ is a
collection of (smooth) vector fields on $M$ which is a $C^\infty(M)$
submodule of $\vectorfields{M}$, the module of vector fields on $M$.
For each $z\in M$, we denote by $\dim_z(\S)$ the dimension of the
subspace of $T_z(M)$ spanned by the values at $z$ of the vector fields
in $\S$.  Now $\dim_z(\S)$ need not be constant.  However, if vector
fields $X_a$ are linearly independent at $z$ they are linearly
independent in a neighbourhood of $z$, which means that
$\dim_{z'}(\S)\geq\dim_z(\S)$ for all $z'$ in a neighbourhood of $z$.
Moreover, $\dim_z(\S)$ has a maximal value on $M$, which we call the
maximal dimension of $\S$, and the set of points at which the maximal
dimension is attained is an open subset of $M$.

A vector field system $\S$ on $M$ restricts to a vector field system
$\S|_U$ on any open subset $U$ of $M$, considered as a submanifold.
On the other hand, any vector field in $\S_U$ may be extended to a
vector field in $\S$ by multiplying it by a bump function whose
support is contained in $U$ (taking advantage of the fact that we are
working in the $C^\infty$ category). So it is permissible to discuss
local aspects of vector field systems in coordinates.

An alternative definition of the term distribution is that a
distribution is a vector field system $\S$ for which $\dim_z(\S)$ is
constant, or for which the maximal dimension is attained everywhere.

We assume that $\V$, mentioned in the opening sentence of the section,
is a distribution in the strict sense. Now suppose that we have a
vector field $F$ on $M$ not belonging to $\V$. We denote by
$\V+[F,\V]$ the collection of vector fields on $M$ which may be
written in the form $V_1+[F,V_2]$ with $V_1,V_2 \in \V$. This is a vector
field system, essentially because for any $f \in C^\infty(M)$ and $V
\in \V$, $f[F,V]=[F,fV] \pmod\V$. We will be concerned with this
vector field system, for a given involutive distribution $\V$ and
vector field $F$, throughout this paper.

\begin{prop}\label{propo1}
Suppose that the vector field $F$ is such that $[F,\V]\cap\V=\{0\}$,
that is, if $V\in\V$ and $[F,V]\in\V$ then $V=0$.  Then the maximal
dimension of $\V+[F,\V]$ is $2n$, and the open subset of $M$ on which
it is attained is dense, that is, its closure is $M$.
\end{prop}

\begin{proof}
We denote the vector field system $\V+[F,\V]$ by $\W$ for convenience.
Clearly the maximal dimension of $\W$ is at most $2n$, and the set of
points $z$ where $\dim_z(\W)=2n$ is open, though it may be empty.
Suppose that $z$ is a point of $M$ with $\dim_z(\W)<2n$.  We show that
there can be no open neighbourhood of $z$ such that $\dim_{z'}(\W)<2n$
for all $z'$ in the neighbourhood.

Since $\V$ is involutive there is a coordinate neighbourhood $U$ of
$z$ and coordinates $(q^a,y^i)$ with $a=1,2,\ldots,m-n$ and
$i=1,2,\ldots,n$ such that the coordinate fields $\partial/\partial
y^i$ span $\V|_U$.  We may set
\[
F|_U=f^a\vf{q^a}+f^i\vf{y^i}
\]
for some smooth functions $f^a,f^i$ on $U$, so that
\[
\left[F|_U,\vf{y^i}\right]=-\fpd{f^a}{y^i}\vf{q^a} \pmod{\V},
\]
and for any $V\in\V_U$, with $V=V^i\partial/\partial y^i$,
\[
[F|_U,V]=-V^i\fpd{f^a}{y^i}\vf{q^a} \pmod{\V}.
\]
Now since $m-n\geq n$ the rank of the matrix $(\partial f^a/\partial
y^i)$ at any point is at most $n$. If it is $n$ at $z$ then
\[
V^i(z)\fpd{f^a}{y^i}(z)=0\quad\implies\quad V^i(z)=0,
\]
the vector fields $[F|_U,\partial/\partial y^i]$ are linearly
independent at $z$, and $\dim_z(\W)=2n$.  So if $\dim_{z'}(\W)<2n$ for
all $z'$ in a neighbourhood of $z$, which we can take to be a
coordinate neighbourhood as above, then the rank of the matrix
$(\partial f^a/\partial y^i)(z')$ is less than $n$, and we can find
functions $V^i$ on a neighbourhood $U$ of $z$, not all vanishing, such
that
\[
V^i\fpd{f^a}{y^i}=0.
\]
Then the vector field $V=V^i\partial/\partial y^i$ on $U$ satisfies
$[F|_U,V]\in\V$.  So (by multiplying by a suitable bump function) we
can find a vector field $V'$ on $M$, not identically zero, with
$V'\in\V$, such that $[F,V']\in\V$, which is a contradiction.  So
every neighbourhood of a point $z$ where $\dim_z(\W)<2n$ must contain a
point $z'$ where $\dim_{z'}(\W)=2n$.  Thus the set of points $z$ where
$\dim_z(\W)=2n$ is nonempty and open, and its closure is $M$; that is,
the set of points where $\dim_z(\W)=2n$ is an open dense subset of
$M$.
\end{proof}

From now on we will assume that the vector field $F$ does indeed
satisfy the condition of Proposition~\ref{propo1}, and we will
restrict our attention to the open subset where $\dim_z(\W)=2n$, that
is to say, we will effectively assume that $\W$ is a distribution.  We
will make the further assumption that the distribution $\W$ is also
involutive.

We will work locally for the rest of this section, and drop explicit
notational reference to the neighbourhood on which we are working.

In the proof of Proposition~\ref{propo1} we showed that if $\{V_i\}$
is a local basis of $\V$ consisting of coordinate fields
$\partial/\partial y^i$ of a local coordinate system $(q^a,y^i)$ and
we set $W_i=[F,V_i]$ then $\{V_i,W_i\}$ is a local basis for $\W$.
Indeed, this will be true for any local basis $\{V_i\}$ of $\V$.  If
we change basis to $\tilde{V}_i=A_i^jV_j$ (where the $A_i^j$ are
locally defined smooth functions and $(A_i^j)$ is nonsingular) then
$W_i$ changes to $\tilde{W}_i=A_i^jW_j+F(A_i^j)V_j$.  If the basis
$\{V_i\}$ is such that $[V_i,V_j]=0$, so that the $V_i$ are coordinate
fields, the necessary and sufficient condition for $\{\tilde{V}_i\}$
also to satisfy $[\tilde{V}_i,\tilde{V}_j]=0$ is that
$A_i^lV_l(A^k_j)=A_j^lV_l(A^k_i)$.

Our next aim is to show that, under the assumptions stated earlier,
one can choose a commuting basis $\{V_i\}$ for $\V$ such that it and
the corresponding $W_i$ satisfy $[V_i,W_j]\in\V$.  For any basis
$\{V_i\}$ we can write
\[
[V_i,W_j]=\alpha^k_{ij}V_k+\beta^k_{ij}W_k.
\]
Notice that if $[V_i,V_j]=0$ then both coefficients are symmetric in
their lower indices:
\[
0=[F,[V_i,V_j]]=[W_i,V_j]+[V_i,W_j]=[V_i,W_j]-[V_j,W_i].
\]
If we change the basis of $\V$ to $\tilde{V}_i=A_i^jV_j$ we have
\begin{align*}
[\tilde{V}_i,\tilde{W}_j]&=[A_i^kV_k,A_j^lW_l+F(A^l_j)V_l]\\
&=(A_i^lV_l(A_j^k)+A_i^lA_j^m\beta^k_{lm})W_k \pmod\V.
\end{align*}
So to make $[\tilde{V}_i,\tilde{W}_j]\in\V$ we want to choose $A^j_i$
such that $V_l(A_j^k)+A_j^m\beta^k_{lm}=0$.  Note that since
$\beta^k_{ij}$ is symmetric, we will then have
$A_i^lV_l(A_j^k)=A_j^lV_l(A_i^k)$, and if the $V_i$ pairwise commute
then the $\tilde{V}_i$ will also pairwise commute.

The equations
\[
\fpd{A_j^k}{y^l}+A_j^m\beta^k_{lm}=0
\]
are linear first-order partial differential equations for the
unknowns $A^i_j$, and admit solutions if and only if their
integrability conditions, which are
\[
\fpd{\beta_{jk}^l}{y^i}-\fpd{\beta_{ik}^l}{y^j}
+\beta_{im}^l\beta_{jk}^m-\beta_{jm}^l\beta_{ik}^m=0,
\]
are satisfied. Now
\[
0=[[V_i,V_j],W_k]=[[V_i,W_k],V_j]+[V_i,[V_j,W_k]],
\]
and
\[
[V_i,[V_j,W_k]]=(V_i(\beta_{jk}^l)+\beta_{im}^l\beta_{jk}^m)W_l \pmod\V.
\]
So it follows from the identity $[V_i,[V_j,W_k]]-[V_j,[V_i,W_k]]=0$
that the integrability conditions are indeed satisfied.  If we take a
solution $A^i_j$ for which the matrix $(A^i_j)$ is nonsingular on a
local cross-section of the $\V$ foliation, for example by taking
$A^i_j=\delta^i_j$ there, then $(A^i_j)$ will be nonsingular on an
open subset containing the cross-section.  We have shown the
following.

\begin{prop}\label{beta}
If both $\V$ and $\W$ are involutive, there is a commuting basis
$\{V_i\}$ of $\V$ such that for all $i,j$, $[V_i,W_j]\in\V$ (where
$W_i=[F,V_i]$).
\end{prop}

The remaining freedom in the choice of commuting basis (such that
$[V_i,W_j]\in\V$ still holds) is to take $A^j_i$ to satisfy
$V_k(A_i^j)=0$.

The condition $[V_i,W_j]\in\V$ says that $W_j$ is invariant under the
action of $\V$, modulo $\V$.

Let us take a coordinate neighbourhood $U$ in $M$, with coordinates
$(q^a,y^i)$ such that $V_i=\partial/\partial y^i$; we may suppose
without essential loss of generality that $U$ is the image of a
product of open subsets $O\subset\R^{m-n}$ and $P\subset\R^n$, where
$0\in P$.  Then $y^i=0$ is a submanifold of $U$ of codimension $n$,
say $N$, and $U$ is fibered over $N$ with fibres the integral
submanifolds of $\V$.  Denote by $\pi:U\to N$ the corresponding
projection.  Then the restriction of $W_j$ to $U$ is projectable to
$N$:\ that is to say, there is a well-defined vector field $\bar{W}_j$
on $N$ which is $\pi$-related to $W_j$.  More generally, a vector
field $X\in\W$, say $X=X^iW_i\pmod{\V}$, is projectable if, and only
if, the coefficients $X^i$ satisfy $V(X^i)=0$ for all $V\in\V$, or
indeed if $V_j(X^i)=0$.  Let us denote by $\bar{\W}$ the distribution
on $N$ spanned by the $\bar{W}_i$, in other words, the distribution
consisting of the projections of projectable vector fields in $\W$.
Then $\bar{\W}$ is involutive, since it $\pi$-related to the
involutive distribution $\W$.  We may therefore choose coordinates
$(t^p,x^i)$ on $N$, where $p=1,2,\ldots,m-2n$, such that the integral
submanifolds of $\bar{\W}$ are given by $t^p=\mbox{constant}$.  Then
with respect to the coordinates $(t^p,x^i,y^i)$ on $U$ we have
\[
V_i=\vf{y^i},\quad
W_i=W^j_i(x)\vf{x^j}\pmod{\V},
\]
where the coefficients $W^j_i$ are everywhere the components of a
nonsingular matrix.  We still have at our disposal the freedom to
change the original basis to $\tilde{V}_i=A_i^j(t,x)V_j$.  If we do so
with $A_i^kW_k^j=\delta_i^j$ then
\[
\tilde{V}_i=A_i^j\vf{y^j},\quad
\tilde{W}_i=\vf{x^i}\pmod{\V}.
\]
If we make a further change of coordinates to
\[
\tilde{t}^p=t^p,\quad\tilde{x}^i=-x^i,\quad
\tilde{y}^i=W^i_j(x)y^j,
\]
then
\[
\vf{\tilde{t}^p}=\vf{t^p}\pmod{\V},\quad
\vf{\tilde{x}^i}=-\vf{x^i}\pmod{\V},\quad
\vf{\tilde{y}^i}=A_i^j\vf{y^j}
\]
We have proved the following result.

\begin{prop}\label{coords}
With $\V$, $\W$, $V_i$, $W_i$ as above, we can find local coordinates
$(t^p,x^i,y^i)$ on $M$, $p=1,2,\ldots,m-2n$, $i=1,2,\ldots,n$ such
that
\[
V_i=\vf{y^i},\quad
W_i=-\vf{x^i}\pmod{\V}.
\]
\end{prop}

We move on now to investigate the form of $F$.  There are two cases to
consider, depending on whether $F$ does or does not belong to $\W$.

\begin{thm} \label{sode} Assume both $\V$ and $\W$ are involutive.

\begin{enumerate}
\item Suppose that $F\in\W$, and assume that the set $N\subset
M=\{z\in M:F(z)\in\V\}$ is nonempty.  Then we may choose coordinates
with respect to which
\[
F=y^i\vf{x^i}+F^i(t,x,y)\vf{y^i};
\]
that is to say, $F$ takes the form of a second-order differential
equation field in terms of the coordinates $(x^i,y^i)$, with the $t^p$
merely behaving as parameters.
\item Suppose that $F$ is everywhere independent of $\W$ (so that in
particular $M>2n$) and that $[F,\W]\subset\W$.  Then we may choose
coordinates with respect to which
\[
F=\vf{t^1}+y^i\vf{x^i}+F^i(t,x,y)\vf{y^i};
\]
that is to say, $F$ takes the form of a time-dependent second-order
differential equation field in terms of the coordinates
$(t^1,x^i,y^i)$, with the $t^p$ with $p>1$ merely behaving as
parameters.
\end{enumerate}

\end{thm}

\begin{proof}

1.\ For the first case, set $F=a^iV_i+b^iW_i$ with respect to a frame
with $[V_i,V_j]=0$ and $[V_i,W_j]\in\V$.  Then
\[
W_i=[F,V_i]=-V_i(a^j)V_j+b^j[W_j,V_i]-V_i(b^j)W_j,
\]
so we must have $V_i(b^j)=-\delta^j_i$.  Now $N$ is the zero level set
of $(b^i)$, and the rank of the Jacobian of the map $M\to\R^n:z\mapsto
(b^i(z))$ is $n$, or in other words the 1-forms $db^i$ are
independent, since $(V_i(b^j))$ is nonsingular.  So $N$ is an immersed
submanifold of $M$ of codimension $n$, and $\V$ is transverse to it.
We may choose coordinates $(t^p,x^i,y^i)$ as in
Proposition~\ref{coords}, such that $V_i=\partial/\partial y^i$, $N$
is given by $y^i=0$, $(t^p,x^i)$ are coordinates on $N$, and
\[
W_i=-\vf{x^i}\pmod{\V}.
\]
With respect to such coordinates set
\[
F=f^i(t,x,y)\vf{x^i}+F^i(t,x,y)\vf{y^i}.
\]
Then from its definition
\[
W_i=-\fpd{f^j}{y^i}\vf{x^j}\pmod{\V},
\]
and therefore
\[
\fpd{f^j}{y^i}=\delta^j_i.
\]
Taking into account the fact that $f^i(t,x,0)=0$ we have
$f^i(t,x,y)=y^i$.

2.\ For the second case, take coordinates as in
Proposition~\ref{coords}, and suppose that
\[
F=\varphi^p(t,x,y)\vf{t^p}\pmod{\W}.
\]
From the assumption that $[F,\W]\subset\W$ it follows that in fact
$\varphi^p$ depends only on the $t^q$.  By assumption the $\varphi^p$
cannot vanish simultaneously, and so by a transformation of the
coordinates $t^p$ we may take
\[
F=\vf{t^1}+f^i(t,x,y)\vf{x^i}+F^i(t,x,y)\vf{y^i}.
\]
Arguing as above we see that
\[
\fpd{f^j}{y^i}=\delta^j_i.
\]
We may only conclude now that $f^i(t,x,y)=y^i+k^i(t,x)$.  However, a
further coordinate transformation $y^i\mapsto y^i+k^i(t,x)$, with
$t^p$ and $x^i$ unchanged, leads to
\[
\vf{t^p}\mapsto\vf{t^p}\pmod{\V},\quad \vf{x^i}\mapsto\vf{x^i}\pmod{\V},
\]
with $\partial/\partial y^i$ unchanged, and so leads to the required
form for $F$.
\end{proof}

The remaining freedom in transforming the coordinates $(x^i,y^i)$ in
the first case, so as to preserve the form of $F$, is
\[
\tilde{x}^i=\tilde{x}^i(x),\quad \tilde{y}^i=\fpd{\tilde{x}^i}{x^j}y^j.
\]
That is to say, the $y^i$ transform like canonical fibre coordinates
on a tangent bundle. From this point of view it is natural to think
of the coordinates $y^i$ in use before the final transformation
leading to Proposition~\ref{coords} as quasi-velocities.

The remaining freedom in the second case is
\[
\tilde{x}^i=\tilde{x}^i(t^1,x),\quad
\tilde{y}^i=\fpd{\tilde{x}^i}{x^j}y^j+\fpd{x^i}{t^1}.
\]
Here the $y^i$ transform like the jet coordinates of the 1-jet bundle
of a manifold fibred over $\R$.

\section{Induced connections}

The coefficients $\beta^k_{ij}$ used in the proof of
Proposition~\ref{beta} have the appearance of the components of a
connection, and the integrability conditions quoted in the proof have
the form of the vanishing of the curvature of this connection.  We
begin this section by explaining in what sense the $\beta^k_{ij}$ are
indeed the components of a flat symmetric connection.

Let $\V$ be an involutive distribution on any manifold $M$.  For any
vector field $X$ on $M$ denote by $X+\V$ the equivalence class of $X$
modulo $\V$, that is, the collection of vector fields differing from
$X$ by an element of $\V$.  For any $V\in\V$, set
\[
D_V(X+\V)=[V,X]+\V.
\]
This is a well-defined operator on equivalence classes, which is
$\R$-linear in both arguments, and for $f\in C^\infty(M)$ satisfies
\[
D_{fV}(X+\V)=fD_V(X+\V),\quad
D_V(fX+\V)=fD_V(X+\V)+V(f)(X+\V).
\]
That is to say, $D$ has connection-like properties. By the Jacobi
identity, for any $V_1,V_2\in\V$
\[
D_{V_1}D_{V_2}(X+\V)-D_{V_2}D_{V_1}(X+\V)-D_{[V_1,V_2]}(X+\V)=0;
\]
that is to say, if $D$ were a connection it would have zero curvature.

More particularly, let $\W$ be another involutive distribution on $M$,
with $\V\subset\W$; then we may restrict $X$ in the construction above
to lie in $\W$.  The same conclusions hold, mutatis mutandis.  We may
think of $\V$ and $\W$ as vector sub-bundles of $T(M)$, and vector
fields in the distributions as sections of the corresponding bundles
$\V\to M$, $\W\to M$.  If $W$ is a section of $\W\to M$ then we may
think of $W+\V$ as a section of the vector bundle $\W/\V\to M$.  Then
(using the terminology of Lie algebroid theory) $D$
is a $\V$-connection on $\W/\V$.

Now take $\dim\W=2n$, $\dim\V=n$, and suppose there is a type $(1,1)$
tensor field $S$ on $\W$ (that is, a section of the bundle
$\W\otimes\W^*\to M$) with the algebraic properties of an almost
tangent structure (so that\ $\im S=\ker S$), with kernel $\V$.  Then $S$
defines an isomorphism between sections of $\W/\V$ and sections of
$\V$.  So we may define a $\V$-connection on $\V$, say $\nabla$, by
\[
\nabla_{V_1}V_2=S(D_{V_1}(W+\V))\quad
\mbox{for any $W\in\W$ such that $S(W)=V_2$.}
\]
That is, $\nabla_V S(W)=S(D_V(W+\V))=S([V,W])$.  This is well-defined
as a $\V$-connection, and has vanishing curvature.  Since $\nabla$ is
a $\V$-connection on $\V$, it makes sense to talk about its torsion.
But for any $W_1,W_2\in\W$,
\begin{eqnarray*}
\lefteqn{\nabla_{S(W_1)}S(W_2)-\nabla_{S(W_2)}S(W_1)-[S(W_1),S(W_2)]}\\
&=&S[S(W_1),W_2]-S[S(W_2),W_1]-[S(W_1),S(W_2)]\\
&=&-([S(W_1),S(W_2)]-S[S(W_1),W_2]-S[W_1,S(W_2)])\\
&=&-N_S(W_1,W_2).
\end{eqnarray*}
That is to say, the torsion vanishes if and only if the formal
Nijenhuis torsion $N_S$ of $S$ (a type $(2,1)$ $\W$-tensor) vanishes.

If $S$ has vanishing Nijenhuis torsion in this sense, and we restrict
attention to any leaf of the involutive distribution $\V$, we obtain a
flat symmetric connection there.

We now show how to construct such a $\W$-tensor $S$ in the case of
interest.

\begin{prop}\label{prop1}  Assume both $\V$ and $\W$ are involutive.
There is a unique type $(1,1)$ $\W$-tensor field $S$ for which
\[
S(V) = 0 \qquad\mbox{and}\qquad S([F,V]) = -V, \qquad V\in\V;
\]
it satisfies $\ker S=\im S =\V$ and $N_S=0$.
\end{prop}

\begin{proof}
Let $\{V_i\}$ be a basis of $\V$, and set $W_i=[F,V_i]$:\ then
$\{V_i,W_i\}$ is a basis for $\W$.  So it is enough to know how $S$
acts on elements of the form $V$ and $[F,V]$.  The definition above is
consistent:\ if $V \in\V$, then also $fV\in\V$ with $f$ a function on
${M}$, and $S([F,fV]) = S(F(f)V) + S(f[F,V]) = fS([F,V]) = -fV$.  We
have $S^2=0$, $\ker S=\im S =\V$.  The formal Nijenhuis torsion $N_S$
obviously vanishes for two elements in $\V$.  Moreover, for any
$V_1,V_2\in\V$,
\[
N_S(V_1,[F,V_2]) = - S[V_1,V_2] = 0
\]
because of the assumed integrability of $\V$.  Likewise, by making use
of the Jacobi identity (and because $[V_1,V_2]\in\V$),
\begin{align*}
N_S([F,V_1],[F,V_2]) &= [V_1,V_2] +S[V_1,[F,V_2]] + S [[F,V_1],V_2]\\
&= [V_1,V_2] + S [F,[V_1,V_2]] = 0.\qedhere
\end{align*}
\end{proof}

In the case of interest, where $S$ is as defined above, the
$\beta^k_{ij}$ are the connection coefficients of this connection with
respect to the basis $\{V_i\}$.

Suppose we have a further distribution $\H$ on $M$, of dimension $n$,
contained in $\W$, and everywhere transverse to $\V$; in other words a
complement to $\V$ in $\W$.  We call such a distribution horizontal.
Then the restriction of $S$ to $\H$ is a $C^\infty(M)$-isomorphism
$\H\to\V$.  For any $V\in\V$, denote by $\hlift{V}$ the unique element
of $\H$ such that $S(\hlift{V})=V$.

We can extend the $\V$-connection $\nabla$ on $\V$ to a
$\W$-connection on $\V$ as follows:\ for any $W\in\W$ and any $V\in\V$
set
\[
\nabla_WV=P_\V([P_\H(W),V])+S([P_\V(W),\hlift{V}]),
\]
where $P_\H$ and $P_\V$ are the projectors on $\H$ and $\V$,
respectively.  The right-hand side belongs to $\V$ and depends
$\R$-linearly on the arguments.  For $f\in C^\infty(M)$
\begin{align*}
\nabla_{fW}V&=P_\V([fP_\H(W),V])+S([fP_\V(W),\hlift{V}])\\
&=f\nabla_WV-V(f)P_V(P_\H(W))-\hlift{V}(f)S(P_\V(W))\\
&=f\nabla_WV
\end{align*}
while
\begin{align*}
\nabla_W(fV)&=P_\V([P_\H(W),fV])+S([P_\V(W),f\hlift{V}])\\
&=f\nabla_WV+P_\H(W)(f)P_\V(V)+P_\V(W)(f)S(\hlift{V})\\
&=f\nabla_WV+(P_\H(W)+P_\V(W))(f)V\\
&=f\nabla_WV+W(f)V.
\end{align*}
So $\nabla$ is a covariant derivative.

If $W\in\V$, say $W=V_1$, then the new definition gives
$\nabla_{V_1}V_2=S([V_1,\hlift{V_2}])$.  According to the old definition,
$\nabla_{V_1}V_2=S([V_1,W])$ for any $W$ such that $S(W)=V_2$.  But
$W=\hlift{V_2}$ is such that $S(W)=V_2$; so the two definitions agree in
this case.  On the other hand, suppose that $W\in\H$ and that $W$ is
projectable in the sense that $[W,\V]\subset\V$ (the horizontal
projection of any projectable vector field is projectable, and the
$W_i$ are projectable as we pointed out before).  Then
$\nabla_WV=[W,V]$.

Assuming as before that $[F,\W]\subset\W$ (which is automatically the
case if $F\in\W$, and is an assumption in Part 2 of Theorem~\ref{sode}
if not), it is possible to define a Lie derivative by $F$ of
$\W$-tensors:\ for example, in the case of a type $(1,1)$ $\W$-tensor
$T$ as the commutator of operators $\ad F$ and $T$:
\[
(\lie{F}T)(W)=[F,T(W)]-T[F,W].
\]

\begin{prop} If $[F,\W]\subset\W$, the vector field $F$ defines a
complement $\H$ of $\V$ in $\W$.
\end{prop}
\begin{proof}
We show now that, with the above definition, $\lie{F} S$ defines two
projection operators on $\W$.  We first show that $(\lie{F}S)^2 =
\id$.  We have, for $V\in \V$,
\[
(\lie{F} S) (V) = [F,S(V)] - S[F,V] = V
\]
and thus $(\lie F S)^2(V)=V$. Also,
\[
(\lie F S) ([F,V]) = [F,S[F,V]] - S[F,[F,V]] = - [F,V] -S[F,[F,V]].
\]
Since $S[F,[F,V]]\in\V$, we have $((\lie{F} S) (S[F,[F,V]]) =
S[F,[F,V]] $, and therefore
\[
(\lie{F} S)^2 ([F,V]) = - (\lie{F} S) ([F,V]) - S[F,[F,V]] = [F,V].
\]
The conclusion is that $P_\H = \onehalf (\id - \lie{F} S )$ and $P_\V =
\onehalf (\id + \lie{F} S )$ are complementary projection operators,
with e.g.\ $P_\V(V)=V$ and $P_\H(V)=0$;  $\H=\im P_\H$ is therefore a
complement to $\V$ in $\W$.
\end{proof}

We will use this complement from now on.

In the case where $M$ is a tangent manifold $T(Q)$, we can take $\V$
to be the canonical vertical distribution, and in particular
$\W=\vectorfields{T(Q)}$.  The connection with covariant derivative
$\nabla$ is then the Berwald connection associated to a system of
autonomous second-order differential equations, see e.g.\ \cite{Berw}
(taking into account the fact that the current connection is expressed
in terms of vertical vector fields rather than vector fields along the
tangent bundle projection).  A similar construction exists for the
case where $M$ is the first jet manifold of a bundle $E\to\R$, and
where the second-order dynamics are time-dependent, see e.g.\
\cite{CMS,MS}.

The Berwald connection can be used to describe special classes of
second-order differential equation fields, such as the ones of
quadratic type we had encountered in the introduction.

The `mixed curvature' of the $\V$-connection $\nabla$ is the (1,2)
$\V$-tensor field $\theta$ given by
\[
\theta(V_1,V_2)V_3=
\nabla_{\hlift{V}_1}\nabla_{V_2}V_3 - \nabla_{V_2}\nabla_{\hlift{V}_1} V_3
- \nabla_{[\hlift{V}_1,V_2]}V_3.
\]

\begin{prop}
Let $\V$ and $\W$ both be involutive and $[F,\W]\subset\W$.  The
necessary and sufficient condition for the existence of coordinates
in which $F$ takes the form of a quadratic second-order differential
equation field is that $\theta=0$.
\end{prop}

\begin{proof}
Let $\hat\V$ denote the set of $V\in \V$ for which the corresponding
$\hlift V$ is projectable, i.e.\ $V$ satisfies $[\hlift V, V_2] \in
\V$ for all $V_2\in\V$.  This set defines a module over the
projectable functions on $M$ (those functions $f$ for which
$V_1(f)=0$, for all $V_1\in\V$).  Alternatively, $V_2\in\hat\V$ if and
only if $\nabla_{V_1} V_2 = 0$, for all $V_1\in\V$.

Let $V_1,V_2\in\hat\V$. Then
\[
  \nabla_{V_3} \nabla_{\hlift{V}_1}V_2=
  -\theta(V_1,V_3)V_2 - \nabla_{[\hlift{V}_1,V_3]}V_2
+ \nabla_{\hlift{V}_1}\nabla_{V_3}V_2= - \theta (V_1,V_3)V_2,
\]
meaning that $\nabla_{\hlift{V}_1}V_2$ is again projectable if and
only if $\theta=0$.

Let $\theta=0$. If we set,
for $V_1,V_2 \in \hat V$,
\[
D_{V_1} V_2 = \nabla_{\hlift{V}_1} V_2,
\]
one easily verifies that the operator $D$ satisfies connection-like
properties with respect to the multiplication of elements of the
module $\hat\V$ with projectable functions $f$:
\[
D_{fV_1} V_2 = f D_{V_1} V_2 \qquad \mbox{and} \qquad
D_{V_1} f V_2 = f D_{V_1} V_2 + \hlift{V}_1(f) V_2.
\]
Remark that $\hlift{V}_1(f)$ is again projectable, since for any
$V_2\in\V$, $V_2 \hlift{V}_1(f) = \hlift{V}_1V_2(f) + [V_2,
\hlift{V}_1](f)=0$.

In the coordinates as defined in Theorem~\ref{sode} (regardless of
whether $F$ lies in $\W$ or not) the connection coefficients are
given by
\[
D_{ \frac{\partial}{\partial y^i} } \fpd{}{y^j} =
\conn kij \fpd{}{y^k}, \qquad
\mbox{where\,\,}  \Gamma^i_j = -\onehalf \fpd{f^i}{y^j} \,
\mbox{\,and\,} \,  \conn ijk = \fpd{\Gamma^i_j}{y^k}.
\]
It is clear from this expression that $\conn kij = \conn kji$, or,
equivalently, that the connection $D$ is symmetric, in the sense that
the torsion
\[
D_{V_1} V_2 - D_{V_2} V_1 - S[\hlift{V}_1,\hlift{V}_2]
\]
vanishes.

We can therefore conclude that the functions $\Gamma^k_{ij}$ are
projectable if and only if $\theta=0$.  For that to be the case, $f^k$
must be of the form $f^k = \conn kij y^iy^j + P^k_i y^i + Q^k$, for
some projectable functions $P^k_i(t,x)$ and $Q^k(t,x)$.
\end{proof}

The advantage of the current description is that the criterion
$\theta=0$ can be verified in any given set of coordinates on $M$.

\section{Global properties}\label{globalsection}

In this section we will address the global bundle structure of a
manifold $M$ in the context of a given involutive distribution and a
vector field $F$, assuming from the start that the set where the
dimension of $\V+[F,\V]$ is maximal is the whole of $M$, or in other
words that $\V+[F,\V]$ is actually a distribution.

Let $M$ be an $m$-dimensional manifold with an involutive (and thus
integrable) $n$-dimensional distribution $\V$.  The foliation of the
distribution defines an equivalence relation on $M$ by declaring two
elements of $M$ to be equivalent if they lie in the same leaf of $\V$.
The quotient of $M$ by means of this equivalence relation, say $Q$,
will have the structure of a differentiable manifold if for every leaf
one can find a smooth embedded local submanifold $N$ through a
point of the leaf, of dimension $m-n$, which has the property that
each other leaf it meets is intersected in only one point.  Then
$\pi_1: M\to Q$ defines a fibration, for which the fibres have dimension $n$,
and for which the distribution $\V$ coincides with the tangents to the
fibres.

We will assume that this condition is satisfied.

In the case of interest, $Q$ comes also equipped with an integrable
distribution.  Indeed, a projectable vector field on $M$ is
$\pi_1$-related to a vector field on $Q$.  Those projectable vector
fields that happen to lie in $\W$ define therefore a distribution on
$Q$, say $\bar \W$, which is involutive by construction.  As above,
the corresponding equivalence relation therefore defines a new
quotient, $T$, again assumed to be a manifold.  We will denote the
corresponding fibration by $\pi_2: Q\to T$.  Alternatively, we could
have defined a fibration by quotienting out the distribution $\W$
from the beginning.  This structure will coincide with the composition
projection $\pi_2 \circ\pi_1: M\to T$.

An almost tangent structure on an even dimensional manifold is a
(1,1)-tensor field $S$ on that manifold, for which the kernel of $S$
at each point coincides with its image.  The almost tangent structure
is said to be integrable if its Nijenhuis torsion vanishes.  If that
is the case then the kernel of $S$ is an involutive distribution.  We
now recall a result from \cite{CT}.  Suppose that the kernel of an
integrable almost tangent structure on a manifold defines a fibration
over some base manifold (as above, by taking the quotient of that
distribution).  Suppose that each fibre is connected and simply
connected, and that there exists a flat connection on each fibre, for
which the fibre is geodesically complete.  Then the manifold is the
total manifold of an affine bundle modelled over the tangent bundle of
the base manifold.

This theorem can be applied to the current setting, if we take a
particular leaf $L_{\W}$ of $\W$ to be the even dimensional manifold
of interest.  The projection $\pi_1$ will project this leaf $L_\W$ of
$\W$ onto a corresponding leaf $L_{\bar\W}$ of $\bar\W$.  Therefore we
may consider the fibration given by the restriction
$\pi_1|_{L_\W}: L_\W \to L_{\bar\W} $.  Vectors that are tangent to
its fibres can be identified with vectors in $\V$, and the fibres
themselves can be identified with leaves of $\V$.  We have defined a
$\W$-tensor field $S$ on $\W$ in Proposition~\ref{prop1}.  It
restricts naturally to an almost tangent structure on $L_\W$ (i.e.\
$S^2=0$, $N_S=0$, and $\ker S=\im S= \V_{L_\W}$).
Moreover, we have seen that the restriction of the $\V$-connection
$\nabla_V$ ($V\in\V$) to a leaf of $\V$ gives a flat connection on
that leaf.  We can conclude therefore:

\begin{thm} \label{thm2}
Suppose that each leaf of $\V$ is connected and simply connected, and
assume that each leaf of $\V$ is geodesically complete with respect to
the restriction of $\nabla_V$ to that leaf.  Then, for each $L_\W$ of
$\W$, $\pi_1|_{L_\W}: L_\W \to L_{\bar\W}$ is an affine bundle,
modelled over the tangent bundle $T(L_{\bar\W})\to L_{\bar\W}$.
Suppose further that the set $N=\{z\in L_\W:F(z)\in\V\}$ is a global
cross-section of $\pi_1|_{L_\W}$:\ then $L_\W$ may be identified with
$T(L_{\bar\W})$ and $N$ with the zero section.
\end{thm}

\begin{cor} \label{cor1}
In case $F\in\W$, and under the assumptions of the previous theorem,
the restriction of $F$ to a certain leaf $L_\W\equiv T(L_{\bar\W})$ will
be a second-order differential equation field on $T(L_{\bar\W})$.
\end{cor}

\begin{proof} This follows easily from the coordinate expression of
$F$.  Restricting $F$ to a leaf is the same as fixing the parameters
$t^p$ to some constant values.
\end{proof}

For completeness, we mention that one may find an alternative
formulation of the theorem of \cite{CT} in \cite{DeF}, where the
global conditions on an (assumed given) symmetric connection are
replaced by global conditions on an (assumed given) vector field.  In
our current framework, the restriction of the vector field $S(F)$ to
$L_\W$ plays the role of that vector field.

In case $F$ does not belong to $\W$, but leaves it invariant, $F$
defines a vector field $\tilde F$ on $T$.  This vector field defines a
1-dimensional involutive distribution on $T$, leading as before to a
fibration $T\to T_0$.  If we assume that the vector field $\tilde
F$ is complete, an integral curve of $\tilde F$ will define a
1-dimensional submanifold $T_1$ of $T$.  In turn, the preimage of
$T_1$ under $\pi_2$ is a collection $E$ of leaves of $\bar\W$ lying
over that integral curve.  We can think of the restriction of $\pi_2$
to $E$ as defining a fibration $\pi_2|_E: E\to T_1$.  Let's denote
its 1-jet bundle by $J^1(E)\to E$.

The distribution $\W_F = \langle F \rangle \oplus \W$ is also
involutive.  Its leaves $L_{\W_F}$ are $(2n+1)$-dimensional manifolds
that are projected by means of $\pi_1$ onto one of the above described
manifolds $E$, corresponding to a certain integral curve (with image
$T_1$) of $\tilde F$.  The fibres of $\pi_1|_{L_{\W_F}}: L_{\W_F} \to
E$ can again be identified with $\V$.  Recall that we had defined a
symmetric flat connection on each leaf of $\V$.

By setting $S(F)=0$ we can extend $S$ to a $\W_F$-tensor field, which
has the property that its restriction to a leaf of $\W_F$ satisfies
$S^2=0$, $N_S = 0$ and $\rank S = n$.  These properties are exactly
those that define, in the terminology of \cite{MP2}, an `almost jet
structure'.  The above mentioned theorem in \cite{CT} has been
generalized to $(2n+1)$-dimensional manifolds with almost jet
structures in \cite{MP2}, where a 1-jet bundle replaces the role
played by a tangent bundle (see the Theorem on page
90 of \cite{MP2}).  We are in the situation that we can apply this
theorem, since each $L_{\W_F}$ is a $(2n+1)$-dimensional manifold with
a $2n$-dimensional distribution $\W|_{L_{\W_F}}$ which is completely
integrable and which is such that $\im S = \V|_{L_{\W_F}} \subset
\W|_{L_{\W_F}}$.  We may therefore conclude:

\begin{thm}
Suppose that each leaf of $\V$ is connected and simply connected, and
assume that each leaf of $\V$ is geodesically complete with respect to
the restriction of $\nabla_V$ to that leaf.  In case $F\notin\W$ and
$[F,\W]\subset\W$, each leaf $L_{\W_F}$ of $\W_F$ is diffeomorphic to the
1-jet bundle $J^1(E)$ of $E\to T_1$.
\end{thm}

\begin{cor}
Under the assumptions of the previous theorem, the restriction of $F$
to a certain leaf $L_{\W_F}\equiv J^1(E)$ will be a time-dependent
second-order differential equation field on $J^1(E)$.
\end{cor}

\begin{proof}
This follows again from the coordinate expression of $F$.  Restricting
$F$ to a leaf is the same as fixing all parameters $t^p$ to some
constant values, except for $t^1$.
\end{proof}

The cases of most obvious interest are those in which the dimension of
$M$ is either $2n$ or $2n+1$ ($n$ being the dimension of $\V$).  We
end this section with a statement of our global results in these
cases, in a form which collects together the assumptions we have made.

\begin{thm}
Let $\V$ be an involutive distribution of dimension $n$ on a manifold
$M$ of dimension $m$, $m=2n$ or $2n+1$; and $F$ a vector field such
that $\V\cap[F,\V]=\{0\}$. Assume that
\begin{itemize}
\item $M$ is fibred over a manifold $Q$ and the leaves of $\V$ are
the fibres of this fibration;
\item each leaf of $\V$ is connected and simply connected;
\item each leaf of $\V$ is geodesically complete with respect to the
flat symmetric connection induced on it (as described in Section~3);
\item $\V+[F,\V]$ is a distribution (necessarily of dimension $2n$).
\end{itemize}

In the case $m=2n$, assume further that
\begin{itemize}
\item the set $\{z\in M:F(z)\in\V\}$ is a global cross-section of
$M\to Q$.
\end{itemize}
Then $M$ may be identified with $T(Q)$ and $F$ with a second-order
differential equation field on $T(Q)$.

In the case $m=2n+1$ assume further that
\begin{itemize}
\item $\W=\V+[F,\V]$ is involutive;
\item $F\notin\W$, $[F,\W]\subset\W$;
\item $F$ is complete;
\item $Q$ is fibred over $\R$, and the leaves of $\bar{\W}$ (the
projection of $\W$ to $Q$) are the fibres of this fibration.
\end{itemize}
Then $M$ may be identified with $J^1(Q)$ and $F$ with a time-dependent
second-order differential equation field on $J^1(Q)$.

\end{thm}

\section{An illustrative example}

The cases in which $\dim M$ is $2n$ or $2n+1$ may be of most obvious
interest, but they are by no means the only cases of interest, as we
now show by an example.

Let $\Q$ be the configuration space of a Lagrangian system with
regular Lagrangian $L$ and assume that $L$ is invariant under the
(free and proper) action of a symmetry Lie group $G$.  In that case,
the Euler-Lagrange field $\Gamma\in\vectorfields{T(\Q)}$ is
$G$-invariant and it can be reduced to a vector field $\check\Gamma$
on $T(\Q)/G$.  The corresponding equations for finding the integral
curves of $\check\Gamma$ are known in the literature as the
`Lagrange-Poincar\'e equations', see e.g.\ \cite{Cendra}.  The
equations determining the reduced vector field $\check\Gamma
\in\vectorfields{T(\Q)/G}$ can be cast in terms of the reduced
Lagrangian $l$ on $T(\Q)/G$.

We will follow here closely the description we have given in
\cite{invlag}.  With the aid of a principal connection on $\Q\to\Q/G$
one may decompose $T(\Q)/G$ into $T(\Q/G) \oplus (\Q\times\g)/G$,
where the action of $G$ on $\g$ is the adjoint action.  In what
follows $(x^i,v^i,w^p)$ are local coordinates on $T(\Q)/G$, where the
$(x^i)$ are coordinates on $\Q/G$, and $(v^i,w^p)$ are fibre
coordinates corresponding to the decomposition.

We will assume that the symmetry group is Abelian.  This has the
advantage that the adjoint action is trivial.  The vector field
$\check\Gamma$ can then be determined from:
\begin{eqnarray*}
&& \check\Gamma({x}^i) = v^i\\
&&\check{\Gamma}\left( \fpd{l}{v^i}\right) -\fpd{l}{x^i} =
K^p_{ik}v^k \fpd{l}{w^p}\nonumber\\
&&\check{\Gamma}\left( \fpd{l}{w^p}\right) =0.
\end{eqnarray*}
Here $K^p_{ik}$ are the components of the curvature of the principal
connection with respect to an invariant basis.  The coordinate
expression of the reduced field is therefore of the form
$\check\Gamma=v^i\partial/\partial x^i + \Gamma^i \partial/\partial
v^i + \Gamma^p \partial/\partial w^p$, where $\Gamma^i$ and $\Gamma^p$
are functions on $T(\Q)/G$.

In the assumption that the matrix $(\partial^2 l/\partial w^p\partial
w^q)$ is everywhere non-singular, and that the relation $\partial
l/\partial w^p=\mu_p$ can therefore be rewritten in the form
$w^p=\rho^p(x,v,\mu)$, we can perform a coordinate transformation
$(x^i,v^i,w^p) \to ({\bar x}^i = x^i,{\bar v}^i = v^i,\mu_p = \partial
l/\partial w^p)$.  The last equation is then simply
$\check{\Gamma}\left( \mu_p\right) =0$, that is to say: the
coordinates $\mu_p$ can be regarded as parameters.  In the new
coordinates the reduced vector field becomes $\check\Gamma={\bar
v}^i\partial/\partial {\bar x}^i + \Gamma^i \partial/\partial {\bar
v}^i + 0 \partial/\partial \mu_p $.

The first two equations determine a system of second-order ordinary
differential equations in the variables $x^i$ with parameters $\mu_p$.
By introducing Routh's (reduced) function
\[
{\mathcal R}^\mu (x,v)= l(x,v, \rho(x,v,\mu)) - \mu_p \rho^p(x,v,\mu)
\]
these equations, restricted to fixed values for the $\mu$\,s, can
equivalently be rewritten as
\[
\check\Gamma\left( \fpd{{\mathcal R}^\mu}{v^i}\right)
-\fpd{{\mathcal R}^\mu}{x^i} = K^p_{ik}v^k\mu_p.
\]
This equation is known as Routh's (reduced) equation for an Abelian
symmetry group, see e.g.\ \cite{MaSch}.

We show that the situation described above is in agreement with the
statements of Theorem~\ref{sode} and Corollary~\ref{cor1}.  Recall
first the definition of the momentum map $J: T(\Q) \to {\g}^*$, where
$\langle J(v),\xi \rangle = (\vlift{\tilde\xi}L)(v)$ (for each
$\xi\in\g$, $\tilde\xi\in\vectorfields{\Q}$ is the corresponding
fundamental vector field).  It is well-known that the map $ T(\Q) \to
\Q\times {\g}^*, (q,v) \mapsto (q,J(v))$ is $G$-equivariant, where the
action of $G$ on $\g^*$ is the coadjoint action.  Therefore, it
reduces to a map ${\check J}: T(\Q)/G \to (\Q \times {\g}^*)/G$.  But
for an Abelian group the adjoint action is trivial, so the coadjoint
action is also trivial.  It follows that $J$ is invariant, and that
the image of the reduced momentum $\check J$ is $\Q/G\times \g^*$.  In
the current coordinates $\check J$ is simply $(\bar x,\bar v,\mu)
\mapsto (\bar x,\mu)$.

Let $M=T(\Q)/G$ and $F=\check\Gamma$.  The distribution $\V = \ker
T\check J$ is clearly involutive.  It has the commuting basis given by
the vector fields $\partial / \partial {\bar v}^i$.  It is easy to see
that $[F,\V]\cap \V = \{0\}$.  The distribution $\W = \V + [F,\V]$ is
spanned by $\{ \partial /\partial {\bar x}^i, \partial /\partial {\bar
v}^i \}$ and is involutive as well.  It is the distribution formed by
the vector fields on $T(\Q)/G$ which are tangent to the level sets of
momentum.  A leaf of $\W$ is thus a particular level set, $\mu_p
=\mu_p^0$.  The corresponding $N$ in the statement of
Theorem~\ref{sode} can be identified with $\im \check J =
\Q/G\times\g^*$, and is non-empty.  So Theorem~\ref{sode} applies.

The quotient space $Q$ of Section~\ref{globalsection} can be identified with $\Q/G \times
\g^*$.  It is trivially fibred over $T=\g^*$; $\bar \W$ is the
distribution formed by the (projections of the) $\partial/\partial
{\bar x}^i$, and any leaf of $\bar\W$ can therefore be identified with
$\Q/G$. According to Theorem~\ref{thm2}, therefore,
the restriction of $F$ to a level set of momentum $\mu=\mu^0$ is a
second-order differential equation field on $\Q/G$ (a leaf of
$\bar\W$); it is of course the one which satisfies Routh's reduced
equation with $\mu_p=\mu^0_p$.

\subsubsection*{Acknowledgements}
The first author is a Postdoctoral Fellow of the Research Foundation
-- Flanders (FWO).  The second author is a Guest Professor at Ghent
University:\ he is grateful to the Department of Mathematics for its
hospitality.  This work is part of the {\sc irses} project {\sc
geomech} (nr.\ 246981) within the 7th European Community Framework
Programme.  We are indebted to W.\ Sarlet for many useful discussions.

\end{document}